\documentclass[a4paper]{article}

\usepackage{url}
\usepackage{color}
\usepackage{enumerate}
\usepackage[mathscr]{euscript}		
\usepackage{dsfont}

\usepackage{psfrag}			
\usepackage{hyperref}

\usepackage{times}
\usepackage{amsmath}
\usepackage{amssymb}
\usepackage{amsthm}

\newcommand{\ind}{\mathds}

\newcommand{\E}{\ensuremath{\mathbb{E}}}
\renewcommand{\P}{\ensuremath{\mathbb{P}}}

\theoremstyle{plain}
\newtheorem{thm}{Theorem}[section]
\newtheorem{corollary}[thm]{Corollary}
\newtheorem{prop}[thm]{Proposition}
\newtheorem{lem}[thm]{Lemma}
\newtheorem{conj}[thm]{Open problem}

\theoremstyle{definition}

\theoremstyle{remark}
\newtheorem{remark}[thm]{Remark}

\newcommand{\N}{\ensuremath{\mathbb{N}}}                       
\newcommand{\Z}{\ensuremath{\mathbb{Z}}}                       
\newcommand{\R}{\ensuremath{\mathbb{R}}}                       

\DeclareMathOperator{\essinf}{ess\,inf}

\numberwithin{equation}{section}

\title{ASYMPTOTIC DIRECTIONS IN RANDOM WALKS IN RANDOM ENVIRONMENT REVISITED}
\author{A. Drewitz $^{1,2,}$\thanks{
Partially supported by the International Research Training Group
``Stochastic Models of Complex Processes''.
}
$^,$\thanks{Partially supported by Iniciativa Cient\'\i fica Milenio P-04-069-F.}\;
\and A.F. Ram\'\i rez $^{1,\dagger,}$\thanks{
Partially supported by Fondo Nacional de Desarrollo Cient\'\i fico
y Tecnol\'ogico grant 1060738.}
}

\date{\today}

\begin{document}
\maketitle
{\footnotesize
\thanks{$^1$ Facultad de Matem\'{a}ticas, Pontificia Universidad Cat\'{o}lica de Chile, Vicu\~{n}a Mackenna 4860, Macul, Santiago, Chile;
{\it e-mail:} \href{mailto:adrewitz@uc.cl}{adrewitz@uc.cl}, \href{mailto:aramirez@mat.puc.cl}{aramirez@mat.puc.cl}}

\thanks{$^2$ Institut f\"ur Mathematik, Technische Universit\"at Berlin, Sekr. MA 7-5, Str. des 17. Juni 136, 10623 Berlin, Germany;
{\it e-mail:} \href{mailto:drewitz@math.tu-berlin.de}{drewitz@math.tu-berlin.de}}
}
\begin{abstract}
 
Recently  Simenhaus in \cite{Si-07}
  proved that  for any  elliptic random walk in random environment,
 transience in the
  neighborhood of a given direction is equivalent to the a.s. existence of a
  deterministic asymptotic direction and to transience in any direction in the
  open half space defined by this asymptotic direction.
 Here we prove an improved   version  of this result and review some open problems.
\end{abstract}
%
{\footnotesize

}

\section{Introduction} \label{introduction}
For
each site $x\in\Z^d$, consider the vector $\omega(x):=\{\omega(x,e): e \in \Z^d, |e|=1\}$ such that
$\omega(x,e)\in (0,1)$ and $\sum_{|e|=1} \omega(x,e)=1$. We call the set of possible
values of these vectors $\mathcal P$ and define  the {\it
  environment} $\omega =\{\omega(x):x\in\Z^d\}\in \Omega:= \mathcal P^{\Z^d}$. We define
a random walk on the random environment $\omega$, as a random walk $\{X_n:n\in\N
\}$ with a transition probability from a site $x\in \Z^d$ to a nearest neighbor site
$x+e$ with $|e|=1$ given by $\omega(x,e)$. Let us call $P_{x,\omega}$ the law of this
random walk starting from site $x$ in the environment $\omega$. Let $\P$
be a probability measure on $\Omega$ such that the coordinates $\{\omega(x)\}$ of
$\omega$ are i.i.d. We call $P_{x,\omega}$ the {\it quenched} law of the random walk in
random environment (RWRE), starting from site $x$. Furthermore,  
we define the {\it averaged} (or {\it annealed}) law of the RWRE starting from $x$  by $P_x:=\int_\Omega P_{x,\omega} \,d\P$.
In this note we discuss some aspects of
RWRE related to the a.s. existence of an asymptotic direction in dimension
$d\ge 2$, briefly
reviewing some of the open questions which have been unsolved and proving
 an improved version of a recent
theorem of Simenhaus on the a.s. existence of an asymptotic direction. 

Some very fundamental and natural questions about this model remain
open. Given a vector $l\in\R^d \backslash \{0\}$, define the event

$$
A_l:=\{\lim_{n\to\infty}X_n\cdot l=\infty\}.
$$
Whenever $A_l$ occurs, we say that the random walk is transient in the
direction $l$. Let also

$$
B_l:=\left\{\liminf_{n\to\infty}\frac{X_n \cdot l}{n}>0\right\}.
$$
Whenever $B_l$ occurs, we say that the random walk is ballistic in the
direction $l$.  We have the following open problem.

\begin{conj}
\label{conj1} In dimensions $d\ge 2$, transience in the direction $l$ implies
ballisticity in the direction $l$. 
\end{conj}
\noindent
Some partial progress related to this
question has been achieved by Sznitman and Zerner \cite{SzZe-99}, and later by
Sznitman in \cite{Sz-00,Sz-01,Sz-02},
which we will discuss below. Under a uniform ellipticity assumption, i.e.
$
\P ( \essinf \min_{\vert e \vert} \omega(0,e) > 0) = 1,
$
the following lemma, which we call Kalikow's zero-one law, was proved by Sznitman and Zerner (cf.
Lemma 1 in \cite{SzZe-99})
using
regeneration times.
Later  Zerner and Merkl \cite{ZeMe-01} derived the corresponding result under the assumption of ellipticity only,
i.e.
$
\P (\min_{\vert e \vert} \omega (0,e) > 0) = 1;
$
cf.
Proposition 3 in \cite{ZeMe-01}.
\begin{lem}[Sznitman-Zerner] For $l\in \R^d\backslash \{0\}$,

$$
P_0(A_{l}\cup A_{-l})=0\quad {\rm or}\quad 1.
$$
\end{lem}
\noindent
On the other hand, in dimension $d=1$ a zero-one law holds, i.e. $P_0(A_l) \in
\{0,1\}$. Zerner and Merkl,
proved the following (see Theorem 1 in \cite{ZeMe-01}).

\begin{thm}[Zerner-Merkl] 
\label{zm}
In dimension $d=2$, for $l\in \R^2\backslash \{0\}$,
$$
P_0(A_{l})=0\quad {\rm or}\quad 1.
$$
\end{thm}
\noindent 
Nevertheless, we still have the following open problem.

\begin{conj}
\label{conj2}
In dimensions $d\ge 3$, for $l\in \R^d\backslash \{0\}$,
$$
P_0(A_{l})=0\quad {\rm or}\quad 1.
$$
\end{conj}
\noindent
Combining Kalikow's zero-one law with the law of large numbers
result of Sznitman and Zerner in \cite{SzZe-99}, Zerner \cite{Ze-02} proved the following theorem.

\begin{thm}[Sznitman-Zerner]
\label{speeds} In dimensions $d\ge 2$, there exists a direction $\nu\in \mathbb{S}^{d-1}$,
$\nu\ne 0$, and $v_1,v_2\in [0,1]$ such that $P_0$-a.s.

$$
\lim_{n\to\infty}\frac{X_n}{n}=v_1\nu \ind{1}_{A_{\nu}}-v_2 \nu \ind{1}_{A_{-\nu}}.
$$

\end{thm}
\noindent 
Indeed, Theorem 3.2.2 of \cite{TaZe-04}, the proof of which can be performed in the same manner
with the assumption of ellipticity only instead of uniform ellipticity,
implies that for $e \in \Z^d$ with $\vert e \vert = 1$ and
\begin{equation} \label{unionTransOne}
P_0(A_e \cup A_{-e}) = 1
\end{equation}
there exist $v_{e}, v_{-e} \in [0,1]$ such that $P_0$-a.s.
\begin{equation} \label{directionalLLN}
\lim_{n \to \infty} \frac{X_n e}{n} = v_e \ind{1}_{A_e} - v_{-e} \ind{1}_{A_{-e}}.
\end{equation}
Combining this with Theorem 1 of \cite{Ze-02} we may omit assumption (\ref{unionTransOne}) and still obtain
(\ref{directionalLLN}). Having (\ref{directionalLLN}) for the elements $e_1, \dots, e_d$ of the standard basis of $\R^d,$
we obtain that $\lim_{n \to \infty} X_n/n$ exists $P_0$-a.s. and may take values in a set of cardinality $2^d.$
Employing the same argument as Goergen in p. 1112 of \cite{Go-06} we now obtain that
$P_0$-a.s. $\lim_{n \to \infty}X_n/n$ takes two values at most. This yields Theorem \ref{speeds}.

Whenever $\lim_{n\to\infty}X_n/|X_n|$ exists $P_0$-a.s.  we call this limit
the {\it asymptotic direction} and we say that  a.s. an asymptotic direction exists. The
existence of an asymptotic direction can already be established assuming some
of the conditions introduced by Sznitman which imply ballisticity.  Let
$\gamma\in (0,1)$ and $l\in \mathbb{S}^{d-1}$. The condition $(T)_\gamma$ holds
  relative to $l$ if for all $l'\in \mathbb{S}^{d-1}$ in a neighborhood of $l$,

\begin{equation}
\label{tg}
\limsup_{L\to\infty}L^{-\gamma}\log P_0(\{X_{T_{U_{l',b,L}}}\cdot l'<0\})<0,
\end{equation}
for all $b>0$, where $U_{l',b,L}=\{x\in\Z^d:-bL<x\cdot l'<L\}$ is a slab and
$T_{U_{l',b,L}}=\inf\{n\ge 0: X_n\notin U_{l',b,L}\}$ is the first exit time
of this slab. On the other hand, one says that condition $(T')$ holds relative
to $l$ if condition $(T)_\gamma$ holds relative to $l$ for every $\gamma \in
(0,1)$. It is known that for each $\gamma\in (0,1)$ condition $(T)_\gamma$
relative to $l$ implies transience in the direction $l$ and that a.s. an
 asymptotic direction exists which is deterministic. Also, for each
$\gamma\in (1/2,1)$, condition $(T)_\gamma$ relative to $l$ implies condition
$(T')$, which in turn implies ballisticity (see \cite{Sz-02}). One of the open problems related
 to condition $(T)_\gamma$ is the following.

\begin{conj}
\label{conj3} 
If (\ref{tg}) is satisfied for $l' \in {\mathbb S}^{d-1},$ then $(T)_\gamma$ holds relative
to $l'.$
\end{conj}

Recently in \cite{Si-07}, Simenhaus established the following theorem which
gives equivalent conditions for the
existence of an a.s. asymptotic direction and showing that transience in a
neighborhood of a
given direction implies that a.s. an asymptotic direction exists.
\begin{thm}[Simenhaus] \label{thmSimenhaus}
The following are equivalent:
\begin{enumerate}
\item
There exists a non-empty open set $O \subset \R^d$ such that 
\begin{equation} \label{condSimenhaus}
P_0(A_l) = 1 \quad \forall \, l \in O.
\end{equation}
\item
There exists $\nu \in \mathbb{S}^{d-1}$ such that $P_0$-a.s.
$$
\lim_{n \to \infty} \frac{X_n}{\vert X_n\vert} = \nu.
$$
\item
There exists $\nu \in \mathbb{S}^{d-1}$ such that $P_0(A_l) = 1$ for all $l \in \R^d$ with $l\cdot \nu > 0.$
\end{enumerate}
\end{thm}
\noindent 
It is natural to wonder if there exists a statement analogous to
Theorem \ref{speeds}, but related only to the existence of a possibly non-deterministic
asymptotic direction. Here we answer affirmatively this question proving the
following generalization of Theorem \ref{thmSimenhaus}.
\begin{thm} \label{TFAE}
The following are equivalent:
\begin{enumerate}
\item
There exists a non-empty open set $O \subset \R^d$ such that
$$
P_0(A_l \cup A_{-l}) = 1 \quad \forall \, l \in O.
$$

\item
There exist $d$ linearly independent vectors $l_1, \dots, l_d \in \R^d$ such that
\begin{equation} \label{transienceAss}
P_0(A_{l_k} \cup A_{-l_k}) = 1 \quad \forall \, k \in \{1, \dots, d\}.
\end{equation}

\item
There exists $\nu \in \mathbb{S}^{d-1}$ with $P_0(A_\nu \cup A_{-\nu}) = 1$ such that $P_0$-a.s.
\begin{equation} \label{twoAsDirections}
\lim_{n \to \infty} \frac{X_n}{\vert X_n \vert} = \ind{1}_{A_{\nu}} \nu - \ind{1}_{A_{-\nu}} \nu.
\end{equation}

\item
There exists $\nu \in \mathbb{S}^{d-1}$ such that
$$
P_0(A_l \cup A_{-l}) = 1
$$
if and only if $l\in \R^d$ is such that $l \cdot \nu \ne 0.$ In this case, $P_0(A_l \Delta A_{\nu}) = 0$ and $P_0(A_{-l} \Delta A_{-\nu}) = 0$ for all  $l$ such that $l\cdot\nu>0$.
\end{enumerate}
\end{thm}
\noindent It should be noted that Theorems \ref{thmSimenhaus} and \ref{TFAE} are interesting only in the
case in which the statement of the Open Problem \ref{conj1} is not proven to be true.
Furthermore, if condition (\ref{transienceAss}) is fulfilled but (\ref{condSimenhaus}) is not, then if asymptotic directions exist
we have to expect at least (and as it turns out at most, see also Proposition 1 in \cite{Si-07})
two of them. However, it is not
known whether condition (\ref{transienceAss}) can be fulfilled while (\ref{condSimenhaus}) is not. In fact, if 
the statement of the Open Problem \ref{conj2}
holds, then the two conditions are equivalent.
Note that due to Kalikow's zero-one law, condition $(d)$ of Theorem \ref{TFAE} yields a complete characterisation
of $P_0(A_l \cup A_{-l})$ for all $l \in \R^d.$ As a consequence of this result, we obtain an a priori sharper version
of $(c)$ in Theorem \ref{thmSimenhaus}:
$$
(c')\; \text{ \it There exists }\nu \in \mathbb{S}^{d-1} \text{ \it such that }
P_0(A_l) = 1 \text{ \it for all } l \in \R^d \text{ \it with } l\cdot \nu > 0 \text{ \it and } P_0(A_l) = 0 \text{ \it if }
l \cdot \nu \leq 0.
$$
This observation and  Theorem \ref{zm} imply that in dimension $d=2$ there are at most three possibilites for the values of the set of probabilities $\{P_0(A_l):l\in\mathbb S^{d-1}\}$: (1) for all $l$, $P_0(A_l)=0$; (2) there exists a $\nu\in\mathbb S^{d-1}$ such that $P_0(A_\nu)=1$ while $P_0(A_l)=0$ for $l\ne\nu$; (3) there exists a $\nu\in\mathbb S^{d-1}$ such that $P_0(A_l)=1$ for $l$ such that $l\cdot\nu>0$ while $P_0(A_l)=0$ for $l$ such that $l\cdot \nu\le 0$.
The following corollary, which can be deduced from Theorem \ref{TFAE}, shows that knowing that there is an $l^*$ such that $P_0(A_{l^*})=1$ and $P_0(A_l)>0$ for all $l$ in a neighborhood of $l^*$, determines the value of $P_0(A_l)$ for all directions $l$.

\begin{corollary} \label{zeroOne}
The following are equivalent:
\begin{enumerate}
\item
There exists $l^* \in \R^d$ and some neighborhood ${\cal U}(l^*)$
such that $P_0(A_{l^*}) = 1$ and $P_0(A_l) > 0$ for all $l \in {\cal U}(l^*).$
\item
There exists $\nu \in \R^d$ such that $P_0(A_l) =1$  for $l$  such that $l \cdot \nu > 0$, while $P_0(A_l)=0$ for $l$ such that $l\cdot\nu\le 0$.
\end{enumerate}
\end{corollary}
\noindent In particular, this shows that in Theorem \ref{thmSimenhaus}, condition $(a)$ can be replaced
by the a priori weaker condition $(a)$ of this corollary. 

In the rest of this paper we prove Theorem \ref{TFAE} and Corollary \ref{zeroOne}. In Section \ref{prelim}
we prove some preliminary results needed for the proofs and in Section \ref{proofs} we apply them to prove the theorem and the corollary.

\section{Preliminary results} \label{prelim}
The implications $(d) \Rightarrow (a) \Rightarrow (b)$ of Theorem \ref{TFAE} are obvious,
so here we introduce the renewal structure and prove some preliminary results
needed to show that $(b) \Rightarrow (c) \Rightarrow (d).$
For $l \in \R^d$ set
$$D_l := \inf\{n \in \N: X_n \cdot l < X_0 \cdot l\}$$
 and for $B \subset \R^d$ define the first-exit time
$$D_B := \inf\{n \in \N: X_n \notin B\};$$
as usual, we set $\inf \emptyset := \infty.$
We also define for $l \in \R^d$ and $s \in [0,\infty)$, 
$$T^l_s := \inf \{ n \in \N : X_n \cdot l > s\}.$$

Due to their linear independence, the vectors $l_1, \dots l_d$ of Theorem \ref{TFAE} (b)
give rise to the following $2^d$ cones:
$$
C_\sigma := \cap_{k=1}^d \{x \in \R^d :  \sigma_k (l_k \cdot x) \geq 0\}, \qquad \sigma \in \{-1, 1\}^d.
$$
Furthermore, for $\lambda \in (0,1]$ and $l \in \R^d \backslash\{0\}$ we will employ the notation
\begin{equation} \label{interpolCone}
C_{\sigma}(\lambda, l) := \cap_{k=1}^d \{x \in \R^d : (\lambda \sigma_k l_k + (1-\lambda) l) \cdot x \geq 0\},
\end{equation}
where the vectors defining the cone are now interpolations of the $\sigma_k l_k$ with $l.$ Note that
$C_{\sigma}(\lambda, l)$ is a non-degenerate cone with base of finite area if and only if the vectors
$\lambda \sigma_k l_k + (1-\lambda)l,$ $k=1, \dots, d,$ are linearly independent. In particular,
$
C_{\sigma}(1, l) = C_{\sigma}
$
for all $\sigma \in \{-1,1\}^d$ and $l.$

We will often choose $\sigma$ such that
$P_0(\cap_{k=1}^d A_{\sigma_k l_k}) > 0,$ which under (\ref{transienceAss}) is possible
since we then have
\begin{equation} \label{sumConesOne}
1 = P_0 (\cap_{k=1}^d A_{l_k} \cup A_{l_{-k}})
= P_0 ( \cup_\sigma \cap_{k=1}^d A_{\sigma_k l_k} )
= \sum_{\sigma} P_0(\cap_{k=1}^d A_{\sigma_k l_k}).
\end{equation}
For a given $\sigma \in \{-1,1\}^d$ which will usually be clear from the context,
we will frequently consider vectors $l \in \R^d$ satisfying the condition
\begin{equation} \label{coneStrictlyContainedDirection}
\inf_{x \in C_\sigma \cap \mathbb{S}^{d-1}} l\cdot x > 0.
\end{equation}
Note here that for
$\sigma$ such that
$P_0(\cap_{k=1}^d A_{\sigma_k l_k}) > 0$ and
$l$ satisfying (\ref{coneStrictlyContainedDirection}),
the inequality
$
P_0(A_l) \geq P_0(\cap_{k=1}^d A_{\sigma_k l_k})
$
implies that the measure
$
P_0(\cdot \vert A_l)
$
is well-defined. For such $l$ we will then show
the existence of a $P_0(\cdot \vert A_l)$-a.s.
asymptotic direction.
The strategy of our proof is based to a significant part on that of Theorem \ref{thmSimenhaus}. 

We start with the following lemma which ensures that if with positive probability the random walk finally
ends up in a cone, then the probability that it does so and never exits a half-space containing this cone
is positive as well.
\begin{lem} \label{coneFinallyHalfspaceAlways}
Let $\sigma \in \{-1,1\}^d$ and $l \in \R^d$ such that (\ref{coneStrictlyContainedDirection}) holds. Then
$$
P_0(\cap_{k=1}^d A_{\sigma_k l_k}) > 0 \Longrightarrow P_0(\cap_{k=1}^d A_{\sigma_k l_k} \cap \{D_l = \infty\}) > 0.
$$
\end{lem}
\begin{proof}
Assume $P_0(\cap_{k=1}^d A_{\sigma_k l_k} \cap \{D_l = \infty\}) =0.$ Then $\P$-a.s.
\begin{equation} \label{contrAss}
P_{0,\omega} (\cap_{k=1}^d A_{\sigma_k l_k} \cap \{D_l = \infty\}) = 0.
\end{equation}
For $y \in \R^d$ with $l\cdot y \geq 0$
this implies 
\begin{equation} \label{zeroProb}
P_{y} (\cap_{k=1}^d A_{\sigma_k l_k} \cap \{D_{\{x : l\cdot x \geq 0\}} = \infty\}) = 0.
\end{equation}
Indeed, if there existed such $y$ with
$\P(\{ \omega \in \Omega \vert P_{y,\omega} (\cap_{k=1}^d A_{\sigma_k l_k} \cap \{D_{\{x : l\cdot x \geq 0\}} = \infty\}) > 0\}) > 0$
then for
$\omega$
such that $P_{y,\omega} (\cap_{k=1}^d A_{\sigma_k l_k} \cap \{D_{\{x : l\cdot x \geq 0\}} = \infty\}) > 0,$
a random walker starting in $0$ would, with positive probability
with respect to $P_{0,\omega},$ hit $y$ before hitting $\{x : l\cdot x < 0\}$  (due to ellipticity) and from there on
finally end up in
$C_\sigma$ without hitting $\{x : l\cdot x < 0\};$ this is a contradiction to (\ref{contrAss}), hence (\ref{zeroProb}) holds.

Choosing a sequence $(y_n) \subset C_\sigma$ such that $l \cdot y_n \to \infty$ as $n \to \infty$ we therefore
get
\begin{align*}
0 &= P_{y_n}(\cap_{k=1}^d A_{\sigma_k l_k} \cap \{ D_{\{x : l\cdot x \geq 0\}} = \infty\})\\
&\geq P_0(\cap_{k=1}^d A_{\sigma_k l_k} \cap \{D_{\{x : l \cdot x \geq -l\cdot y_n\}} = \infty\}) \to P_0(\cap_{k=1}^d A_{\sigma_k l_k})
\end{align*}
as $n \to \infty.$
To obtain the inequality we employed the translation invariance of $\P$ as well as the monotonicity of events.
\end{proof}
The following lemma will be employed to set up a renewal structure; it can in some way be seen
as an analog to Lemma 1 of \cite{Si-07}.

\begin{lem} \label{coneAlwaysLemma}
Let $\sigma \in \{-1,1\}^d$ such that
$P_0(\cap_{k=1}^d A_{\sigma_k l_k}) > 0.$ Then for each $l$ such that 
(\ref{coneStrictlyContainedDirection}) holds, one has
\begin{equation} \label{coneAlways}
P_0(\{D_{C_\sigma(\lambda, l)} = \infty\}) > 0
\end{equation}
for $\lambda > 0$ small enough.
\end{lem}
\begin{proof}
Lemma \ref{coneFinallyHalfspaceAlways}
implies $P_0(\cap_{k=1}^d A_{\sigma_k l_k} \cap \{D_l = \infty\}) > 0.$
Due to the ellipticity of the walk and the independence of the environment we therefore obtain
\begin{equation} \label{alwaysConeEvent}
P_0 \Big(\{X_1 \cdot l > 0\} \bigcap \cap_{k=1}^d A_{\sigma_k l_k}(X_{1+\cdot} -X_1)
\cap \{D_l(X_{1+\cdot} -X_1) = \infty\} \Big) > 0,
\end{equation}
where we name explicitly the path $X_{1+\cdot} -X_1$ to which the corresponding events
$A_{\sigma_k l_k}$ and $D_l$ refer.
Each path of the event in (\ref{alwaysConeEvent})
is fully contained in $C_{\sigma}(\lambda,l)$ for $\lambda > 0$ small enough.
Thus, the continuity from above of $P_0$ yields
\begin{equation} \label{bigIntersectionPosProb}
P_0 \Big(\{D_{C_\sigma(\lambda,l)} = \infty\}
\cap \{X_1 \cdot l > 0\} \bigcap 
\cap_{k=1}^d A_{\sigma_k l_k}(X_{1+\cdot}-X_1) \cap \{D_l ( X_{1+\cdot}-X_1) = \infty\} \Big) > 0
\end{equation}
for all $\lambda > 0$ small enough.
\end{proof}

Employing Lemma \ref{coneAlwaysLemma}, 
for $\sigma \in \{-1,1\}^d$ with $P_0(\cap_{k=1}^d A_{\sigma_k l_k})> 0$
as in \cite{Si-07} we can introduce a cone renewal structure,
where we choose $l \in \R^d$ such that (\ref{coneStrictlyContainedDirection}) is fulfilled and
the cone to work with is $C_l := C_\sigma(\lambda,l),$ where we fixed
$\lambda>0$ small enough as in the statement of Lemma \ref{coneAlwaysLemma}. Note that for fixed $l$ the
set $C_\sigma(\lambda, l)$ is indeed a cone as long as $\lambda > 0$ is chosen small enough
(since the defining vectors in (\ref{interpolCone}) are linearly independent).
We define
$$
S_0^{l}:= T_0^{l},\quad
R_0^{l}:= D_{X_{S_0^l} +C_l} \circ \theta_{S^{l}_0} + S^{l}_0, \quad
M_0^{l}:= \max \{X_n \cdot {l} : 0 \leq n \leq R_0^{l}\}
$$
and inductively for $k \geq 1:$
$$
S_k^{l}:= T^{l}_{M_{k-1}^{l}}, \quad  R_k^{l} := D_{X_{S^l_{k}} +C_l} \circ \theta_{S^{l}_{k-1}} + S^{l}_{k}, \quad
M_k^{l}:= \max \{X_n \cdot {l} : 0 \leq n \leq R_k^{l}\},
$$
where for $x \in \Z^d$ by $x+C_l$ we denote the cone $C_l$ shifted such that its apex lies at $x.$
Furthermore, set
$$
K^{l} := \inf \{ k \in \N : S^{l}_k < \infty, R^{l}_k = \infty\}
$$
as well as
$$
\tau_1^{l} := S^{l}_{K^{l}},
$$
i.e. $\tau_1^l$ is the first time at which the walk reaches a new maximum in
direction $l$ and never exits the cone $C_l$ shifted to $X_{\tau_1^l}.$
We define inductively the sequence of cone renewal times with respect to $C_l$ by
$$
\tau_k^{l} := \tau_1^{l} (X_{\cdot + \tau_{k-1}^{l}} - X_{\tau_{k-1}^{l}}) + \tau_{k-1}^l
$$
for $k \geq 2.$

The following lemma shows that under the conditions of Lemma \ref{coneAlwaysLemma} the sequence
$\tau_k^{l}$ is well-defined on $A_l.$ It can be seen as an analog to Proposition 2 of \cite{Si-07}.
\begin{lem} \label{RTwellDefLem}
Let
$\sigma \in \{-1,1\}^d$ such that $P_0(\cap_{k=1}^d A_{\sigma_k l_k}) > 0$ and choose $l$ and $\lambda$
such that (\ref{coneStrictlyContainedDirection}) and (\ref{coneAlways}) hold.
Then $P_0(\cdot \vert A_l)$-a.s. one has $K^{l} < \infty.$
\end{lem}
\begin{proof}
Employing Lemma \ref{coneAlwaysLemma}, the proof takes advantage of
standard renewal arguments and is analogous to the proof of Proposition 2
in \cite{Si-07} or Proposition 1.2 in \cite{SzZe-99}.
\end{proof}

\begin{lem} \label{distributionLem}
Let
$\sigma \in \{-1,1\}^d$ such that $P_0(\cap_{k=1}^d A_{\sigma_k l_k}) > 0$ and choose $l$
and $\lambda$ such that (\ref{coneStrictlyContainedDirection}) and (\ref{coneAlways}) hold.
Then
$((X_{\tau_1^{l} \wedge \cdot}, \tau_1^{l}), \dots,
(X_{(\tau_k^{l} + \cdot)\wedge \tau^{l}_{k+1}} - X_{\tau_k^{l}}, \tau_{k+1}^{l} - \tau_k^{l})), \dots$
are independent under $P_0(\cdot \vert A_l)$ and for $k \geq 1,$
$((X_{(\tau_k^{l} + \cdot)\wedge \tau_{k+1}^{l}} - X_{\tau_k}^{l}), \tau_{k+1}^{l} - \tau_k^{l})$
under $P_0( \cdot \vert A_l)$ is distributed like
$(X_{\tau_1^{l} \wedge \cdot}, \tau_1^{l})$ under $P_0(\cdot \vert \{ D_{C_l} = \infty\} ).$
\end{lem}
\begin{proof}
The proof is analogous to the proof of Corollary 1.5 in \cite{SzZe-99}.
\end{proof}

The following lemma has been derived  in Simenhaus'
thesis \cite{Si-08}  (Lemma 2 in there). Here we state it and prove it under a
slightly weaker assumption.

\begin{lem} \label{lemma2Mod}
Let
$\sigma \in \{-1,1\}^d$ such that $P_0(\cap_{k=1}^d A_{\sigma_k l_k}) > 0$ and choose $l \in \Z^d$ and $\lambda$
such that (\ref{coneStrictlyContainedDirection}) and (\ref{coneAlways}) hold and the g.c.d. of the coordinates of $l$ is $1.$ Then
$$
E_0(X_{\tau_1^l} \cdot l \vert \{D_{C_l} = \infty\})
= \frac{1}{P_0(\{D_{C_l} = \infty\}\vert A_l) \lim_{i \to \infty} P_0 (\{T^l_{i-1} < \infty, X_{T^l_{i-1}} \cdot l = i\})}
< \infty
$$
and
\begin{equation} \label{expExistence}
E_0 (X_{\tau_1^l} \vert \{ D_{C_l} = \infty\})
\end{equation}
is well-defined.
\end{lem}

\begin{remark}
A fundamental consequence of working with the cone renewal structure
instead of working with slabs is
the existence of (\ref{expExistence}), see Proposition (\ref{LLNinCone}) also.
\end{remark}

\begin{proof}
The proof leans on the proof of Lemma 3.2.5 in \cite{TaZe-04} which is due to Zerner. Due to the
strong Markov property and the independence and translation invariance of the environment we have for $i > 0:$
\begin{align}
\begin{split} \label{splitAtRenewal}
P_0(\{ \exists k \geq 1 : X_{\tau^l_k} \cdot l = i\} \cap A_l)
&= \sum_{x \in \Z^d, l \cdot x =i} \E P_{0,\omega}
(\{T_{i-1}^l < \infty, X_{T^l_{i-1}} = x, D_{C_l + X_{T^l_{i-1}}} \circ \theta_{T^l_{i-1}} = \infty\})\\
&= \sum_{x \in \Z^d, l\cdot x =i} \E P_{0,\omega}
(\{T_{i-1}^l < \infty, X_{T^l_{i-1}} = x\}) P_{x,\omega}(\{D_{C_l + x} = \infty\})\\
&= P_0(\{T_{i-1}^l < \infty, X_{T^l_{i-1}} \cdot l = i\}) P_0(\{D_{C_l}= \infty\}).
\end{split}
\end{align}
At the same time using
$\{\tau^l_1 < \infty\} = A_l,$ a fact which is proven similarly to
Proposition 1.2 of \cite{SzZe-99}, we compute
\begin{align}
\begin{split} \label{splitForRenewal}
\lim_{i \to \infty} &P_0(\{ \exists k \geq 1 : X_{\tau_k^l} \cdot l = i\} \vert A_l)\\
&= \lim_{i \to \infty} P_0(\{ \exists k \geq 2 : X_{\tau_k^l} \cdot l = i\}\vert A_l)\\
&= \lim_{i \to \infty} \sum_{n \geq 1}
P_0(\{ \exists k \geq 2 : X_{\tau^l_k} \cdot l = i\} \cap \{X_{\tau_1^l} \cdot l = n\} \vert A_l)\\
&= \lim_{i \to \infty} \sum_{n \geq 1}
P_0 (\{ \exists k \geq 2 : (X_{\tau^l_k}-X_{\tau^l_1}) \cdot l = i-n\} \cap \{X_{\tau^l_1} \cdot l = n\} \vert A_l)\\
&= \lim_{i \to \infty} \sum_{n \geq 1}
P_0 (\{ \exists k \geq 2 : (X_{\tau^l_k}-X_{\tau^l_1}) \cdot l = i-n\} \vert A_l)
P_0 (\{X_{\tau_1^l} \cdot l = n\} \vert A_l),
\end{split}
\end{align}
where to obtain the last equality we took advantage of Lemma \ref{distributionLem}. Blackwell's renewal theorem in combination with Lemma
\ref{distributionLem} now yields
$$
\lim_{i \to \infty} P_0(\{ \exists k \geq 2 : (X_{\tau^l_k} -X_{\tau_1^l}) \cdot l = i-n\} \vert A_l)
= \frac{1}{E_0(X_{\tau^l_1} \cdot l \vert \{D_{C_l} = \infty\})}
$$
and thus (\ref{splitForRenewal}) implies
$$
\lim_{i \to \infty} P_0(\{ \exists k \geq 1 : X_{\tau_k^l} \cdot l = i\}\vert A_l)
= \frac{1}{E_0(X_{\tau_1^l} \cdot l \vert \{D_{C_l} = \infty\})}.
$$
Therefore, taking into consideration (\ref{splitAtRenewal}) we infer
\begin{equation} \label{expectationEst}
E_0(X_{\tau_1^l} \cdot l \vert \{D_{C_l} = \infty\})
= \frac{1}{P_0(\{D_{C_l} = \infty\}\vert A_l) \lim_{i \to \infty} P_0 (\{T^l_{i-1} < \infty, X_{T^l_{i-1}} \cdot l = i\})}.
\end{equation}
It remains to show that the right hand side of (\ref{expectationEst}) is finite.
Writing $l_{max} := \max \{\vert l_1 \vert, \dots, \vert l_d \vert\}$ for the maximum of the absolute
values of the coordinates
of $l$ we have
$$
\sum_{i = k}^{k + l_{max}-1} P_0(\{T^l_{i-1} < \infty, X_{T^l_{i-1}} \cdot l = i\})
\geq \sum_{i=k}^{k+ l_{max}-1} P_0(\{T^l_{i-1} < \infty, X_{T^l_{k-1}} \cdot l = i\})
\geq P_0(A_l), \quad \forall k \in \N,
$$
where the first inequality follows since $\{X_{T^l_{k-1}} \cdot l = i\} \subseteq \{X_{T^l_{i-1}} \cdot l = i\}$ for
all $k \in \N$ and $i \in \{k, \dots, k+l_{max}-1\}.$
This now yields
$
\lim_{i \to \infty} P_0(\{T^l_{i-1} < \infty, X_{T^l_{i-1}} \cdot l=i\}) \geq l_{max}^{-1} P_0(A_l) > 0,
$
whence due to (\ref{expectationEst}) we obtain
\begin{equation} \label{expFinite}
E_0 (X_{\tau^l_1} \cdot l \vert \{D_{C_l} = \infty\}) < \infty.
\end{equation}
Since on $\{D_{C_l} = \infty\}$ there exists a constant $C > 0$ such that
$\vert X_{\tau^l_1} \vert \leq C X_{\tau_1^l} \cdot l$, we infer as a direct consequence of (\ref{expFinite})
that (\ref{expExistence}) is well-defined.
\end{proof}

We can now employ the above renewal structure to obtain an a.s. constant
asymptotic direction on $A_l.$
\begin{prop} \label{LLNinCone}
Let
$\sigma \in \{-1,1\}^d$ such that $P_0(\cap_{k=1}^d A_{\sigma_k l_k}) > 0$ and choose $l \in \Z^d$ and $\lambda$
such that (\ref{coneStrictlyContainedDirection}) and (\ref{coneAlways}) hold and the g.c.d of the coordinates
of $l$ is $1.$
Then $P_0(\cdot \vert A_l)$-a.s.
$$
\lim_{n \to \infty} \frac{X_n}{\vert X_n \vert} =
\frac{E_0(X_{\tau_1^l} \vert \{D_{C_l} = \infty\})}{\vert E_0(X_{\tau_1^l} \vert \{D_{C_l} = \infty\})\vert}.
$$
\end{prop}
\begin{remark} \label{lIndep}
In particular, this proposition implies that the limit does not depend on the particular choice of
$l$ nor $\lambda$ (for $\lambda$ sufficiently small). Note that the independence of $l$ stems from the
fact that if $l_1, l_2$ satisfy (\ref{coneStrictlyContainedDirection}) we have
$
P_0(A_{l_1} \cap A_{l_2}) > 0.
$
\end{remark}
\begin{proof}
Due to Lemmas \ref{coneAlwaysLemma} to \ref{lemma2Mod} we may apply the law of large numbers to 
the sequence $(X_{\tau_k^l})_{k \in \N}$ yielding
$$
\frac{X_{\tau^l_k}}{k} \to E_0 (X_{\tau^l_1} \vert \{D_{C_l} = \infty\}) \quad
P_0(\cdot \vert A_l)-a.s., \quad k \to \infty,
$$
and hence
$$
\frac{X_{\tau^l_k}}{\vert X_{\tau^l_k} \vert } \to \frac{E_0(X_{\tau^l_1} \vert \{D_{C_l} = \infty\})}
{\vert E_0(X_{\tau^l_1} \vert \{D_{C_l} = \infty\}) \vert}
\quad P_0(\cdot \vert A_l)-a.s., \quad k \to \infty.
$$
Using standard methods to estimate the intermediate terms (cf. p. 9 in \cite{Si-07}) one obtains
$$
\lim_{n \to \infty} \frac{X_n}{\vert X_n \vert}
= \frac{E_0(X_{\tau^l_1} \vert \{D_{C_l} = \infty\})}{\vert E_0(X_{\tau^l_1} \vert \{D_{C_l} = \infty\}) \vert}
\quad P_0(\cdot \vert A_l)-a.s.
$$
\end{proof}
 The following two results will be needed to obtain results about transience in directions
orthogonal to the asymptotic direction.
\begin{lem} \label{liminfLemma}
Let $(Y_n)_{n \in \N}$ be an i.i.d. sequence on some probability space $({\cal X}, {\cal F}, P)$ with
expectation $EY_1=0$ and variance $EY_1^2 \in (0,\infty].$
Then, for $S_n := \sum_{k=1}^n Y_k$ we have $P (\{\liminf_{n \to \infty} S_n = -\infty\}) = 
P (\{\limsup_{n \to \infty} S_n = \infty\}) =1.$
\end{lem}
\begin{proof}
We only prove $P(\{\liminf_{n \to \infty} S_n = -\infty\}) = 1,$ the remaining equality is
proved in an anolog way.
Setting $\varepsilon := (-\essinf Y_1/2) \wedge 1$ one can show
for all $x \in \R,$ using the strong Markov property
at the entrance times of $S_n$ to the interval $[x,x+\varepsilon],$ 
that $P(\{\liminf_{n \to \infty} S_n \in [x, x+\varepsilon]\}) = 0.$ This then implies
$P(\{\liminf_{n \to \infty} S_n = \pm \infty\}) = 1.$
But Kesten's result
in \cite{Ke-75} yields $\liminf_{n \to \infty} S_n/n > 0$
$P(\cdot \cap \{\liminf_{n \to \infty} S_n = \infty\})$-a.s.,
while by the strong law
of large numbers we have $\lim_{n \to \infty} S_n/n = 0$ $P$-a.s. This
yields $P(\{\liminf_{n \to \infty} S_n = \infty\}) = 0$ and hence finishes the proof.
\end{proof}

\begin{lem} \label{asDirPosProbLem}
Let $l \in \R^d$ such that
\begin{equation} \label{asDirPosProb}
P_0(\{\lim_{n \to \infty} X_n/\vert X_n\vert = l\}) > 0.
\end{equation}
Then, for $l^* \in \R^d$ such that $l^* \cdot l = 0$ one has $P_0((A_{l^*} \cup A_{-l^*}) \cap A_l) = 0.$
\end{lem}
\begin{proof}
We choose a basis
$l_1, \dots, l_d$ of $\R^d$ and $\sigma$ such that $l$ is contained in the interior
of the cone $C_\sigma$ corresponding to $l_1, \dots, l_d$ and  (\ref{coneStrictlyContainedDirection}) is satisfied.
Furthermore, by (\ref{asDirPosProb}) and Lemma \ref{coneAlwaysLemma} we may choose $\lambda$ such that
condition (\ref{coneAlways}) is satisfied for the corresponding cone $C_\sigma(\lambda, l).$
Lemma \ref{RTwellDefLem} yields that the sequence
$(\tau_k^l)_{k \in \N}$ is well defined and Lemmas \ref{distributionLem} and \ref{lemma2Mod}
yield that under $P_0(\cdot \vert A_l)$ the sequence
$( (X_{\tau_2^l} - X_{\tau_1^l})\cdot l^*, (X_{\tau_3^l} - X_{\tau_2^l})\cdot l^*, \dots)$
is i.i.d. with expectation $0,$
the latter being due to the validity
of Lemma \ref{lemma2Mod} as well as (\ref{twoAsDirections}) and $l^* \cdot l = 0.$
Indeed, Proposition \ref{LLNinCone} yields
$$
E_0 (X_{\tau_1^l} \cdot l^* \vert \{D_{C_\sigma(\lambda, l)} = \infty \})
= \vert E_0 (X_{\tau_1^l} \vert \{ D_{C_\sigma (\lambda, l)} = \infty\}) \vert
\underbrace{
\lim_{k \to \infty} \frac{X_{\tau_k^l}}{\vert X_{\tau_k^l} \vert}
}_{=l} \cdot
l^* = 0 \quad P_0(\cdot \vert A_l)-a.s.
$$
Applying Lemma \ref{liminfLemma} to the sequence 
$((X_{\tau_2^l} - X_{\tau_1^l})\cdot l^*, (X_{\tau_3^l} - X_{\tau_2^l})\cdot l^*, \dots)$
yields $P_0((A_{l^*} \cup A_{-l^*}) \cap A_l) = 0.$
\end{proof}

\section{Proof of Theorem \ref{TFAE} and Corollary \ref{zeroOne}} \label{proofs}
\subsection{Proof of Theorem \ref{TFAE}}

We first prove that condition (b) implies (c).
Note that due to Lemma \ref{asDirPosProbLem} and (\ref{coneStrictlyContainedDirection}),
we obtain $P_0(\{\lim_{n \to \infty} X_n/\vert X_n \vert \in \cup_{\sigma} \partial C_\sigma\}) = 0.$ 
We now choose $\sigma$ such that $P_0(\cap_{k=1}^d A_{\sigma_k l_k}) > 0$ and
observe $\cup_{l} \{x \in \R^d: l\cdot x > 0\}
= \text{int} \cup_{\sigma^* \not= -\sigma} C_{\sigma^*};$
here, with ``int'' we denote the interior of a set and
the union is taken over all vectors $l \in \Z^d$ that satisfy (\ref{coneStrictlyContainedDirection}) 
and for which the g.c.d. of the coordinates of $l$ is $1.$
Hence,
letting $l$ vary over all such vectors,
Proposition \ref{LLNinCone} yields
$P_0(\cdot \vert \cup_{\sigma^* \not= -\sigma}\cap_{k=1}^d A_{\sigma^*_k l_k})$-a.s.
\begin{equation} \label{LLNlDirection}
\lim_{n \to \infty} \frac{X_n}{\vert X_n \vert} =
\frac{E_0(X_{\tau_1^{l}}
\vert \{D_{C_{l}} = \infty\})}{\vert E_0(X_{\tau_1^{l}} \vert \{D_{C_{l}} = \infty\})\vert} =: \nu,
\end{equation}
which due to Remark \ref{lIndep} is independent of the respective $l$ chosen.
Now if $P_0(\cup_{\sigma^* \not= -\sigma}\cap_{k=1}^d A_{\sigma^*_k l_k}) = 1$
this finishes the proof and the result is equivalent to Theorem \ref{thmSimenhaus}
obtained in \cite{Si-07}.
Thus, assume
\begin{equation} \label{transienceProb}
P_0(\cup_{\sigma^* \not= -\sigma}\cap_{k=1}^d A_{\sigma^*_k l_k}) \in (0,1).
\end{equation}
In the same manner as before we obtain for any
$l' \in \Z^d$ with coordinates of g.c.d. $1$ and satisfying (\ref{coneStrictlyContainedDirection}) with $\sigma$ replaced
by $-\sigma$
\begin{equation} \label{LLNminusLDirection}
\lim_{n \to \infty} \frac{X_n}{\vert X_n \vert} =
\frac{E_0(X_{\tau_1^{l'}}
\vert \{D_{C_{l'}} = \infty\})}{\vert E_0(X_{\tau_1^{l'}} \vert \{D_{C_{l'}} = \infty\})\vert}
\end{equation}
$P_0(\cdot \vert \cap_{k=1}^d A_{-\sigma_k l_k})$-a.s.
with hopefully self-explaining notations.
Now
Proposition 1 of \cite{Si-07} states that if two elements $\nu \not= \nu'$ of $\mathbb{S}^{d-1}$
occur with positive probability each with respect to $P_0$ as asymptotic directions, then $\nu = -\nu'.$
Thus, (\ref{LLNlDirection}) to (\ref{LLNminusLDirection}) imply that the limit in (\ref{LLNminusLDirection})
equals $-\nu,$ and (\ref{LLNminusLDirection}) holds $P_0(\cdot \vert A_{-\nu})$-a.s. This yields $(c).$

Now with respect to the implication $(c) \Rightarrow (d)$ note that the only thing that is not
obvious at a first glance is that $l\cdot \nu = 0$ implies $P_0(A_l \cup A_{-l}) = 0.$ However,
Lemma \ref{asDirPosProbLem} yields $P_0((A_l \cup A_{-l})\cap (A_\nu \cup A_{-\nu})) = 0$ which
due to $P_0(A_\nu \cup A_{-\nu}) = 1$ yields the desired result.
\subsection{Proof of   Corollary \ref{zeroOne}}
We only have to prove $(a) \Rightarrow (b).$
Given $(a),$ Theorem \ref{TFAE} yields the existence of $\nu \in \mathbb{S}^{d-1}$ such that 
\begin{equation} \label{transNuMeasOne}
P_0(A_\nu \cup A_{-\nu}) = 1
\end{equation}
and (\ref{twoAsDirections}) holds.

Now if $l^* \cdot \nu \not= 0$ then
$P_0(A_\nu \cap A_{l^*})=1$ or $P_0(A_{-\nu} \cap A_{l^*})=1,$ respectively, and 
hence $P_0(A_\nu) = 1$ or $P_0(A_{-\nu})=1,$ which due to
Theorem \ref{thmSimenhaus} finishes the proof.
Thus, assume 
\begin{equation} \label{orthogonality}
l^* \cdot \nu = 0
\end{equation}
from now on.
Then Lemma \ref{asDirPosProbLem} yields $P_0((A_{l^*} \cup A_{-l^*}) \cap (A_\nu \cup A_{-\nu})) =0$ which due to
(\ref{transNuMeasOne}) implies $P_0(A_{l^*} \cup A_{-l^*}) = 0,$ a contradiction to assumption $(a).$
\vspace{1 em}

{\bf Acknowledgments.} We  thank Fran\c cois Simenhaus for reading a
preliminary version of this paper and for making several useful comments on it.

\newpage

\def\polhk#1{\setbox0=\hbox{#1}{\ooalign{\hidewidth
  \lower1.5ex\hbox{`}\hidewidth\crcr\unhbox0}}}


\begin{thebibliography}{Goe06}

\bibitem[Goe06]{Go-06}
Laurent Goergen.
\newblock Limit velocity and zero-one laws for diffusions in random
  environment.
\newblock {\em Ann. Appl. Probab.}, 16(3):1086--1123, 2006.

\bibitem[Kes75]{Ke-75}
Harry Kesten.
\newblock {Sums of stationary sequences cannot grow slower than linearly.}
\newblock {\em Proc. Am. Math. Soc.}, 49:205--211, 1975.

\bibitem[Sim07]{Si-07}
Fran\c{c}ois Simenhaus.
\newblock Asymptotic direction for random walks in random environments.
\newblock {\em Ann. Inst. H. Poincar\'{e}}, 43(6):751--761, 2007.

\bibitem[Sim08]{Si-08}
Fran\c{c}ois Simenhaus.
\newblock {\em Marches Al\'{e}atoires en Milieux Al\'{e}atoires -- \'{E}tude de
  quelques Mod\`{e}les Multidimensionnels}.
\newblock PhD thesis, Universit\'{e} Paris 7 - Denis Diderot, 2008.

\bibitem[SZ99]{SzZe-99}
Alain-Sol Sznitman and Martin Zerner.
\newblock A law of large numbers for random walks in random environment.
\newblock {\em Ann. Probab.}, 27(4):1851--1869, 1999.

\bibitem[Szn00]{Sz-00}
Alain-Sol Sznitman.
\newblock {Slowdown estimates and central limit theorem for random walks in
  random environment.}
\newblock {\em J. Eur. Math. Soc. (JEMS)}, 2(2):93--143, 2000.

\bibitem[Szn01]{Sz-01}
Alain-Sol Sznitman.
\newblock {On a class of transient random walks in random environment.}
\newblock {\em Ann. Probab.}, 29(2):724--765, 2001.

\bibitem[Szn02]{Sz-02}
Alain-Sol Sznitman.
\newblock An effective criterion for ballistic behavior of random walks in
  random environment.
\newblock {\em Probab. Theory Related Fields}, 122(4):509--544, 2002.

\bibitem[TZ04]{TaZe-04}
Simon Tavar{\'e} and Ofer Zeitouni.
\newblock {\em Lectures on probability theory and statistics}, volume 1837 of
  {\em Lecture Notes in Mathematics}.
\newblock Springer-Verlag, Berlin, 2004.
\newblock Lectures from the 31st Summer School on Probability Theory held in
  Saint-Flour, July 8--25, 2001, Edited by Jean Picard.

\bibitem[Zer02]{Ze-02}
Martin P.~W. Zerner.
\newblock A non-ballistic law of large numbers for random walks in i.i.d.\
  random environment.
\newblock {\em Electron. Comm. Probab.}, 7:191--197 (electronic), 2002.

\bibitem[ZM01]{ZeMe-01}
Martin P.~W. Zerner and Franz Merkl.
\newblock A zero-one law for planar random walks in random environment.
\newblock {\em Ann. Probab.}, 29(4):1716--1732, 2001.

\end{thebibliography}
\end{document}